\newcommand{\abstand}{\vspace{1em}}
\newcommand{\retour}{\vspace{-0.7em}}

\newtheorem{coro}{Corollary}
\newtheorem{theo}{Theorem}

\newcounter{axiomzahl}
\newcommand{\axiomvorne}{}
\newcommand{\axiomhinten}{.}
\newenvironment{axiom}{\begin{list}{\rm \axiomvorne
\arabic{axiomzahl}\axiomhinten}{\usecounter{axiomzahl}
 \itemsep0em plus 0.1em \parsep0em plus 0.1em \rm}}{\end{list}}

\newcommand{\allgemeinaxiom}[3]{
      \renewcommand{\axiomvorne}{#1}
      \renewcommand{\axiomhinten}{#2}
      \begin{axiom}{#3}\end{axiom}
      \renewcommand{\axiomvorne}{}\renewcommand{\axiomhinten}{.}}

\newcommand{\zitat}[4]{\bibitem{#1}{\sc #2}: {\em #3\/}.
#4.\vspace{-0.7em}}

\newcommand{\PG}[2]{\mbox{$\mbox{{\rm PG}}(#1,#2)$}}
\newcommand{\AG}[2]{\mbox{$\mbox{{\rm AG}}(#1,#2)$}}
\newcommand{\PL}{\mbox{$(\Pcal,\Lcal)$}}
\newcommand{\PLL}{\mbox{$(\Pcal',\Lcal')$}}

\newcommand{\GF}[1]{\mbox{$\mbox{{\rm GF}}(#1)$}}

\newcommand{\abb}[3]{\mbox{$#1\,:\,#2\rightarrow#3$}}
\newcommand{\Abb}[5]{\mbox{$#1\,:\,#2\rightarrow#3,\;#4\mapsto #5$}}
\newcommand{\pabb}[3]{\mbox{$#1\,:\,#2\succ\!\rightarrow#3$}}
\newcommand{\ppfeil}{\succ\!\rightarrow}
\newcommand{\im}{\mbox{\rm im\,}}

\newcommand{\dom}{\mbox{\rm dom\,}}
\newcommand{\ex}{\mbox{\rm ex\,}}

\newcommand{\spn}{\mbox{{\rm span\,}}}

\newcommand{\qu}[1]{\overline{#1}}

\newcommand{\abf}{{\bf a}}
\newcommand{\bbf}{{\bf b}}

\newcommand{\Mbf}{{\bf M}}

\newcommand{\Vbf}{{\bf V}}
\newcommand{\Wbf}{{\bf W}}

\newcommand{\Hcal}{{\cal H}}

\newcommand{\Lcal}{{\cal L}}
\newcommand{\Mcal}{{\cal M}}

\newcommand{\Ocal}{{\cal O}}
\newcommand{\Pcal}{{\cal P}}

\newcommand{\Scal}{{\cal S}}
\newcommand{\Tcal}{{\cal T}}

\newcommand{\Xcal}{{\cal X}}
\newcommand{\Ycal}{{\cal Y}}

\documentclass[a4paper,11pt]{article}
\usepackage{amsfonts}
\title{Weak Linear Mappings -- A Survey}
\author{Hans Havlicek}
\date{}
\begin{document}
\maketitle
\thispagestyle{empty} 
\noindent
There are various concepts of {\em structure preserving mappings}\/ in
geometry. It is the aim of the present paper to give a survey on geometrical
characterizations of some of those mappings. We discuss the results for
projective spaces in some detail and report on generalizations to other
spaces.

We shall come across {\em partially defined mappings}. The notation
$\pabb{\varphi}{\Pcal}{\Pcal'}$ is used in order to point out that $\varphi$
is a mapping with {\em domain}\/ $\dom\varphi\subset\Pcal$ and {\em image
set} $\im\varphi\subset\Pcal'$. The {\em exceptional set}\/ of $\varphi$
(with respect to $\Pcal$) is given as $\ex\varphi:=\Pcal \setminus
\dom\varphi$. The mapping $\varphi$ is said to be {\em globally defined}\/
(with respect to $\Pcal$) provided that $\ex\varphi$ is empty. The notation
$\abb{\varphi}{\Pcal}{\Pcal'}$ is maintained for globally defined mappings
only.

If we are given any subset $\Mcal\subset\Pcal$, then $\{X^\varphi \mid
X\in\Mcal\cap\dom\varphi \}$ is a well--defined set. By abuse of notation, it
is written as $\Mcal^\varphi$. Hence $\Mcal^\varphi=\emptyset$ exactly for
$\Mcal\subset\ex \varphi$.

\section{Projective Spaces}

\subsection{Weak Semilinear Mappings of Vector Spaces}
\label{semilinear}

Throughout this paper $\Vbf$, $\Vbf'$ denote right vector spaces over
(not necessarily commutative) fields $K$, $K'$, respectively.
We generalize the well--known concept of a semilinear mapping as follows: A
{\em weak semilinear mapping}%
   \footnote{Such a mapping should not be mixed up with a {\em generalized
   semilinear mapping}\/ \cite{mach75} which yields a {\em homomorphism}\/ of
   projective spaces.}
\abb{f}{\Vbf}{\Vbf'} with respect to a mapping $\abb{\zeta}{K}{K'}$ is
additive and $\zeta$--homogeneous, i.e.\
   \begin{displaymath}
   (\abf + \bbf)^f = \abf^f + \bbf^f \mbox{ and }
   (\abf x)^f=\abf^fx^\zeta \mbox{ for all }\abf,\bbf\in \Vbf,\,x\in K.
   \end{displaymath}
If $f\neq 0$, then $\zeta$ is a monomorphism which is uniquely determined by
$f$.

Given such an $f\neq 0$ there exists a right vector space $\Wbf$ over
$F:=K^\zeta$ and a semilinear bijection $\abb{g}{\Vbf}{\Wbf}$ with respect to
$\zeta$ (regarded as isomorphism $K\rightarrow F$). By the universal property
of the tensor product $\Wbf\otimes_F K'$, there is a $K'$--linear mapping
$\abb{h}{\Wbf\otimes_F K'}{\Vbf'}$ such that $\abf^f=(\abf^g\otimes 1)^h$ for
all $\abf\in\Vbf$.

The one-- and two--dimensional subspaces of $\Vbf$ are the points and lines of
the projective space on $\Vbf$ which is denoted by
$(\Pcal(\Vbf),\Lcal(\Vbf))$. More generally, we put
$\Pcal(\Mbf):=\{ \abf K \mid \abf\in \Mbf\setminus \{ 0 \} \}$
for $\Mbf\subset \Vbf$; $\Pcal(\Wbf)$ etc.\ is defined likewise.
Each weak semilinear mapping \abb{f}{\Vbf}{\Vbf'} determines a mapping of
points
   \begin{equation}\label{p(f)}
   \pabb{\varphi}{\Pcal(\Vbf)}{\Pcal(\Vbf')},\;
   \abf K\mapsto (\abf^f)K'\mbox{ for all } \abf K \in
   \Pcal(\Vbf\setminus\ker f).
   \end{equation}
The exceptional set of $\varphi$ is the subspace $\Pcal(\ker f)$. It follows
that $\varphi$ has the property
\begin{equation}
   \label{schwachverbindungstreu}
   (\Xcal\vee \Ycal)^\varphi \subset
   \spn(\Xcal^\varphi)\vee\spn(\Ycal^\varphi)
   \mbox{ for all subspaces }\Xcal, \Ycal\subset\Pcal(\Vbf).
   \end{equation}
We remark that $\varphi$--images of subspaces need not be subspaces.

Suppose now that $\zeta$ is bijective, whence $f$ is semilinear. If
$\Xcal\subset\Pcal$ is a subspace, then so is $\Xcal^\varphi$. Moreover,
(\ref{schwachverbindungstreu}) improves to
\begin{eqnarray}
   \label{verbindungstreu}
   (\Xcal\vee \Ycal)^\varphi &=& \Xcal^\varphi\vee\Ycal^\varphi
   \mbox{ for all subspaces }\Xcal, \Ycal\subset\Pcal(\Vbf).
   \end{eqnarray}

\subsection{Definition and Examples of Linear Mappings}

In the sequel \PL\ and \PLL\ denote arbitrary projective spaces. Property
(\ref{verbindungstreu}) is adopted in the following purely geometric
definition:

A {\em linear mapping}\/%
   \footnote{This name has been used in Descriptive Geometry for more than
   100 years. Some authors use other terminologies.}
of projective spaces is a mapping \pabb{\varphi}{\Pcal}{\Pcal'} (defined on
some subset of $\Pcal$) satisfying the following conditions:
\allgemeinaxiom{(L}{)}
   {
   \item ${(\{X\}\vee \{Y\})^\varphi} = \{X\}^\varphi\vee\{Y\}^\varphi$
   for all $X,Y\in \Pcal$, $X\neq Y$.
   \item If $X^\varphi=Y^\varphi$ for distinct points $X,Y\in\dom\varphi$,
then there is at least one exceptional point on the line $\{X\}\vee \{Y\}$.
   }
Condition (L2) is only used to rule out the possibility that {\em all}\/
points of a line are mapped to the same point. There are three
possibilities for the image of a line $l\in\Lcal$:

$\#(l\cap\ex\varphi)\geq 2$: Then $l\subset\ex\varphi$, i.e.\ we have an {\em
exceptional line}.

$\#(l\cap\ex\varphi)=1$: Then all points of $l\cap \dom\varphi$ are mapped to
the same point.

$\#(l\cap\ex\varphi)=0$: Then $\varphi\mid l$ is injective and
$l^\varphi\in\Lcal'$.

\abstand
\noindent
We give some examples of linear mappings:
\vspace{-0.7em}
   \begin{axiom}
   \item The mapping (\ref{p(f)}) is linear, if $\abb{f}{\Vbf}{\Vbf'}$ is
   semilinear.
   \item Let $\Scal\subset\Pcal$ be a subspace. Write $\Pcal/\Scal$ for the
   set of all subspaces of the form $\Scal\vee\{X\}$, where
   $X\in\Pcal\setminus\Scal$. Then $\Pcal/\Scal$ is the set of points of the
   quotient space of \PL\ modulo $\Scal$. The {\em canonical projection}\/
      \begin{displaymath}
      \pabb{\psi_\Scal}{\Pcal}{\Pcal/\Scal},
      \;
      X\mapsto \Scal\vee\{X\} \mbox{ for all } X\in\Pcal\setminus\Scal
   \end{displaymath}is
   a surjective linear mapping with $\ex\psi_\Scal=\Scal$.

   \item With $\Scal$ as above, the {\em canonical injection}\/
   $\Abb{\iota_\Scal}{\Scal}{\Pcal}{X}{X}$
   is a globally defined linear mapping of the projective space on  $\Scal$
   in $\Pcal$.
   \item A {\em projection}\/ is based on two complementary subspaces of
   $\Scal,\Tcal$ of $\PL$ as follows:
      \begin{displaymath}
      \pabb{\pi}{\Pcal}{\Tcal},\; X\mapsto X^\pi
      \mbox{ with }
      \{X^\pi\}:=(\Scal\vee\{X\})\cap\Tcal \mbox{ for all }
      X\in\Pcal\setminus\Scal.
      \end{displaymath}
   This $\pi$ is a surjective linear mapping with $\ex\pi=\Scal$.
   \item Any collineation of projective spaces is a globally defined
   bijective linear mapping.
   \item
   Suppose that \pabb{\varphi}{\Pcal}{\Pcal'} and
   \pabb{\varphi'}{\Pcal''}{\Pcal'''} are linear mappings with
   $\im\varphi\subset\Pcal''$. Then \pabb{\varphi\varphi'}{\Pcal}{\Pcal'''}
   is also a linear mapping (possibly an empty mapping).
   \item
   Let \pabb{\varphi}{\Pcal}{\Pcal'} be linear. The hyperplanes of \PL\ and
   \PLL\ are, respectively, the points of the dual projective spaces
   $(\Pcal^\ast,\Lcal^\ast)$ and $(\Pcal'^\ast,\Lcal'^\ast)$.
   The {\em extended $\varphi$--preimage}\/ of any hyperplane
   $\Hcal'\subset\Pcal'$ is given by $\ex\varphi \cup \Hcal^{\varphi^{-1}}$.
   Define
   \begin{displaymath}
   \pabb{\varphi^T}{\Pcal'^\ast}{\Pcal^\ast},\;
   \Hcal'\mapsto\ex\varphi \cup \Hcal'^{\varphi^{-1}}
   \mbox{ for all }\Hcal'\in\Pcal'^\ast
   \mbox{ with }
   \im\varphi\not\subset\Hcal'
   \end{displaymath}
   This {\em transpose mapping}\/ of $\varphi$ is linear \cite{faur-froe95},
   \cite{havl85}.
   \end{axiom}

\subsection{Brauner's Theorem}

   \begin{theo} {\sc (Brauner \cite{brau73})}
   Let $\pabb{\varphi}{\Pcal}{\Pcal'}$ be a mapping of projective spaces.
   \vspace{-0.7em}
   \begin{enumerate}
   \item
   If $\varphi$ satisfies (L1), then the exceptional subset $\ex\varphi
   \subset \Pcal$ and the image set\/ $\im\varphi \subset \Pcal'$ are
   subspaces.
   \vspace{-0.7em}
   \item
   If $\varphi$ satisfies (L1) and\/ $\#\im \varphi\neq 1$, then $\varphi$ is
   linear.
   \vspace{-0.7em}
   \item
   Each linear mapping $\varphi$ is decomposable into the canonical
   projection $\Pcal\ppfeil\Pcal/\ex\varphi$, a collineation of this quotient
   space onto the subspace\/ $\im\varphi$, and the canonical injection of\/
   $\im\varphi$ in $\Pcal'$.
   \end{enumerate}
   \end{theo}
As a matter of fact {\sc Brauner} did not discuss mappings $\varphi$
satisfying (L1) and $\#\im\varphi=1$. It is immediate that here $\ex\varphi$
can be any subspace of \PL\ other than $\Pcal$.

By the second fundamental theorem of projective geometry we obtain:
   \begin{coro} \label{coro1}
   Each linear mapping \pabb{\varphi}{\Pcal(\Vbf)}{\Pcal(\Vbf')} such that\/
   $\im\varphi$ contains a triangle is induced by a semilinear mapping
   \abb{f}{\Vbf}{\Vbf'} which is determined to within a non--zero factor in
   $K'$.
   \end{coro}
The existence of a triangle in $\im\varphi$ is needed in Corollary
\ref{coro1} to avoid ``degenerate'' cases:
\vspace{-0.7em}
\begin{axiom}
   \item If $\im\varphi$ is a line, then such an $f$ needs not exist, since a
   collineation of $1$--dimensional projective spaces is just a bijection of
   their point sets. However, when $\#K\in\{2,3,4\}$, then $\varphi$ is
   nevertheless induced by a semilinear mapping.
   \item If $\im\varphi$ is a singleton or empty, then $K$ and $K'$ need not
   be isomorphic. If $K$ and $K'$ are assumed to be isomorphic, then
   $\varphi$ can be induced by a semilinear mapping.
\end{axiom}
\vspace{-0.7em}
We remark that under the assumptions of Corollary \ref{coro1} the transpose
of the semilinear mapping $f$ induces the transpose of $\varphi$.

A geometric characterization of linear mappings of real projective spaces is
due to {\sc Rehbock} \cite{rehb26a}, \cite{rehb26b}. {\sc Timmermann}
\cite{timm73} characterizes the projections in projective spaces by
conditions similar to our (L1). \cite[4.5]{brau-I76} and a paper by {\sc
Faure} and {\sc Fr\"o\-licher} \cite{faur-froe93} contain proofs of Brauner's
Theorem.
A characterization of a linear mapping $\varphi$ of a projective space in
its dual space in terms of {\em two}\/ mappings is given by {\sc Faure} and
{\sc Fr\"o\-licher} \cite{faur-froe95}. The two mappings are $\varphi$ and
$\varphi^T|\Pcal$, where $\Pcal$ is identified with a subspace of the bidual
projective space. When these two mappings coincide, then one obtains a
possibly degenerate quasipolarity. See also {\sc Lenz \cite{lenz54},
\cite{lenz57a}}.

\subsection{Definition and Examples of Weak Linear Mappings}

Formula (\ref{schwachverbindungstreu}) motivates the following definition:  A
{\em weak linear mapping}\/ of  projective spaces is a mapping
\pabb{\varphi}{\Pcal}{\Pcal'} (defined on some subset of $\Pcal$) satisfying
\allgemeinaxiom{(WL}{)}
   {
   \item ${(\{X\}\vee \{Y\})^\varphi} \subset \{X\}^\varphi\vee\{Y\}^\varphi$
   for all $X,Y\in \Pcal$, $X\neq Y$.
   }
By (WL1), there are four possibilities for the image of a line $l\in\Lcal$:

   $\#(l\cap\ex\varphi)\geq 2$: Then $l\subset\ex\varphi$, i.e.\ we have
   an {\em exceptional line}.

   $\#(l\cap\ex\varphi)=1$: Then all points of $l\cap \dom\varphi$ are mapped
   to the same point.

   $\#(l\cap\ex\varphi)=0$ and all points of $l$ are mapped to the same
   point.

   $\#(l\cap\ex\varphi)=0$ and $\varphi\mid l$ is injective, but $l^\varphi$
   is not necessarily a line.

\noindent
Thus now, in contrast to a linear mapping, we do not rule out the possibility
that {\em all}\/ points of a line are mapped onto the same point.

Each mapping $\abb{\varphi}{\Pcal}{\Pcal'}$, which is globally defined,
injective, and collinearity--preserving, is a weak linear mapping. If
non--collinearity of points is also preserved, then $\varphi$ is called an
{\em embedding}. An embedding is called {\em strong}, if independent points
in \PL\ always go over to independent points in \PLL.

\abstand

\noindent
We give some examples of weak linear mappings:
\vspace{-0.7em}
   \begin{axiom}
   \item Formula (\ref{p(f)}) yields a weak linear mapping.
   \item Each linear mapping is also a weak linear mapping.
   \item Let $L/K$ be a field extension and let $\Vbf$ be a vector space over
   $K$. Then
   \begin{displaymath}
   \Abb{\varphi}{\Pcal(\Vbf)}{\Pcal(\Vbf\otimes_K L)}
   {\abf K}{(\abf\otimes 1)L}
   \end{displaymath}
   is a strong embedding. This $\varphi$ is a collineation if, and only if,
   $K=L$ or $\dim\Vbf\leq 1$.
   \item The following embedding is not strong ({\sc Brezuleanu} and
   {\sc R\u{a}dulescu} \cite{brez-radu84}): Let $L/K$ be a field extension
   such that there exist elements $1,y_0,y_1,y_2\in L$ which are linearly
   independent in the {\em left}\/ vector space $L$ over $K$.
   Then define
   \begin{displaymath}
   \Abb{\varphi}{\Pcal(K^4)}{\Pcal(L^3)}
      {(a_0 ,\; a_1 ,\; a_2 ,\; a_3 )K}
      {(a_0+y_0 a_3 ,\; a_1+y_1 a_3 ,\; a_2+y_2 a_3)L}.
   \end{displaymath}
   \item Each product of weak linear mappings is a weak linear mapping.
   \end{axiom}

\subsection {A Fundamental Theorem for Weak Linear Mappings}

   \begin{theo}\label{theo-wl}
   {\sc (Faure} {\rm and} {\sc Fr\"olicher \cite{faur-froe94}, Havlicek
   \cite{havl94})} Each weak linear mapping
   $\pabb{\varphi}{\Pcal(\Vbf)}{\Pcal(\Vbf')}$
   such that $\im\varphi$ contains a triangle is induced by a weak semilinear
   mapping $\abb{f}{\Vbf}{\Vbf'}$ which is determined to within a non--zero
   factor in $K'$.
   \end{theo}
We infer from the decomposition of a weak semilinear mapping in
\ref{semilinear} the following
   \begin{coro}\label{coro2}
   Each weak linear mapping $\pabb{\varphi}{\Pcal}{\Pcal'}$ of desarguesian
   projective spaces such that $\im\varphi$ contains a triangle can be
   decomposed into a strong embedding $\Pcal\rightarrow\Pcal''$ in some
   projective space $(\Pcal'',\Lcal'')$ and a linear mapping
   $\Pcal''\ppfeil\Pcal'$.
   \end{coro}
Corollary \ref{coro2} does not remain true if the hypothesis on $\im\varphi$
is dropped:
If $\im\varphi$ does not contain a triangle of \PLL, then $\ex\varphi$ is
still a subspace. For each point $Y\in\im\varphi$ the extended preimage
$\ex\varphi \cup \{Y\}^{\varphi^{-1}}$ is a subspace. Thus
$\Pcal/\ex\varphi$ is partitioned into subspaces. However, those subspaces are
not necessarily isomorphic. Conversely, any subspace $\Scal\subset\Pcal$ and
any partition of $\Pcal/\Scal$ into subspaces determines a weak linear
mapping in some projective line.
   \begin{coro}
   Let $\pabb{\varphi}{\Pcal(\Vbf)}{\Pcal(\Vbf')}$ be given as in
   Theorem \ref{theo-wl}. Furthermore, assume that each monomorphism
   $K\,\rightarrow\,K'$ is surjective. Then $\varphi$ is a linear mapping.
   \end{coro}
For example, each monomorphism of $\Bbb R$ is surjective \cite[p.\
88]{benz92}, whereas $\Bbb C$ admits non--surjective monomorphisms \cite[p.\
114]{benz92}.

A recent result due to {\sc Kreuzer} \cite{kreu97x} says that each pappian
projective space is embeddable in a pappian projective plane. Other results
on embeddings and specific examples are due to {\sc Benz} \cite[p.\
109]{benz92}, {\sc Brezuleanu} and {\sc R\u{a}dulescu} \cite{brez-radu84},
{\sc Brown} \cite{brow88}, {\sc Carter} and {\sc Vogt} \cite{cart-vogt80a},
\cite{cart-vogt80b}, {\sc Havlicek} \cite{havl88}, and {\sc Limbos}
\cite{limb80a}, \cite{limb81}, \cite{limb82}.

Another special class of weak linear mappings are {\em semicollineations}\/
of projective spaces, i.e.\ globally defined bijective mappings which
preserve the collinearity, but not necessarily the non--collinearity of
points. {\sc Ceccherini} \cite{cecc67} has given an example of a
semicollineation of a $4$--dimensional projective space onto a
non--desarguesian projective plane. Cf. also {\sc Bernardi} and {\sc Torre}
\cite{bern-torr84}, and {\sc Maroscia} \cite{maro70}. It seems to be an open
problem, if each semicollineation of desarguesian projective spaces with
dimensions $\geq 2$ is a collineation.

\subsection{Local Characterizations}

The definition of a linear mapping $\varphi$ uses all points of $\Pcal$.
The following ``local'' results characterize a linear mapping only in terms
of its domain or in terms of a subset of its domain.

Firstly,  we state the following variant of axiom (L1):
\allgemeinaxiom{(L}{')}
   {
   \item ${(\{X\}\vee \{Y\})^\varphi} = \{X\}^\varphi\vee\{Y\}^\varphi$
   for all $X,Y\in \dom\varphi$, $X\neq Y$.
   }
   \begin{theo}{\sc (S\"orensen \cite{soer85})}\label{soerensen}
   Let $\pabb{\varphi}{\Pcal}{\Pcal'}$ be a mapping of projective spaces
   satisfying (L1'). Then the following conditions hold true:
   \retour
   \begin{enumerate}
   \item
   The image set $\im\varphi \subset \Pcal'$ is a subspace.
   \retour
   \item
   If $\spn(\ex\varphi)\neq\Pcal$ and if $\#\im\varphi\geq 2$, then the
   exceptional subset $\ex\varphi \subset \Pcal$ is a subspace.
   \retour
   \item
   If $\spn(\ex\varphi)\neq\Pcal$ and if $\im\varphi$ contains a triangle,
   then $\varphi$ is a linear mapping.
   \end{enumerate}
   \end{theo}
The subsequent examples illustrate that the assumptions in Theorem
\ref{soerensen} are essential:
   \retour
   \begin{axiom}
   \item Let $\abb{\varepsilon}{\Pcal'}{\Pcal}$ be a non--surjective strong
   embedding with $\Pcal'\neq\emptyset$, e.g., the canonical injection of a
   proper subspace $\Pcal'\subset\Pcal$. Define
   $\pabb{\varphi}{\Pcal}{\Pcal'}$ as inverse mapping of $\varepsilon$ with
   domain $\Pcal'^\varepsilon$. This $\varphi$ satisfies (L1'), but
   $\ex\varphi\subset\Pcal$ is not a subspace.
   \item Let $\#\Pcal > 1$ and $\#\Pcal'=1$. There exists an $\Mcal\subset
   \Pcal$ other than a subspace. Define $\pabb{\varphi}{\Pcal}{\Pcal'}$ by
   $\ex\varphi=\Mcal$. Hence $\#\im\varphi=1$ and (L1') is trivially true.
   \item Let \PL\ be a projective space containing a line $t_1$ with more than
   three points. Choose a complement $\Scal_1$ of $t_1$ and let $\Scal$ be a
   hyperplane containing $\Scal_1$. Write \pabb{\varphi_1}{\Pcal}{t} for the
   projection onto $t$ with exceptional subspace $\Scal_1$ and put
   $t:=t_1\setminus\Scal$. The mapping
   $\Abb{\varphi}{\Pcal\setminus\Scal}{t}{X}{X^{\varphi_1}}$
   satisfies (L1') with $\PLL=(t,\{t\})$, but clearly (L1) is violated. Thus
   $\varphi$ is not a linear mapping \cite{soer85}.
   \end{axiom}
Secondly, axioms (L1) and (L2) are modified as follows:
\allgemeinaxiom{(L}{'')}
   {
   \item $X,Y,Z$ collinear implies $X^\varphi, Y^\varphi, Z^\varphi$
   collinear, for all $X,Y,Z\in\dom\varphi$.
   \item If $X^\varphi=Y^\varphi$ for distinct points $X,Y\in\dom\varphi$,
   then $\varphi\mid ((\{X\}\vee \{Y\})\cap\dom\varphi)$ is a constant
   mapping.
   }
Moreover, we need some topological tools: A projective space \PL\ is said to
carry a {\em linear topology}, if each line $x\in\Lcal$ is endowed with a
non--trivial%
   \footnote{We rule out the coarsest and the finest topology.}
topology $T_x$ such that all perspectivities between intersecting lines are
continuous. A subset $\Ocal$ of $\Pcal$ is called {\em linearly open}, if
$\Ocal\cap x$ is open in the topological space $(x,T_x)$ for all lines
$x\in\Lcal$ \cite{bial74}.
If \PL\ is a topological projective space \cite[Ch.\ 23]{buek95}, then the
induced topologies on the lines yield a linear topology. Each open set
$\Mcal\subset \Pcal$ is also linearly open. However, there are linear
topologies that do not arise in this way. An example is given by the {\em
cofinite topology}\/ in a projective space with infinite order: A subset $m$
of a line $a$ is defined to be open, if $a\setminus m$ is finite.

   \begin{theo}\label{frank}
   {\sc (Frank \cite{frank92})} Let \PL\ and \PLL\ be projective spaces
   satisfying the minor axiom of Desargues. Suppose that
   $\pabb{\varphi}{\Pcal}{\Pcal'}$ is a mapping satisfying (L1'') and (L2'')
   such that $\im\varphi$ contains a triangle. Each of the following
   conditions is sufficient for the existence of a unique linear mapping
   \pabb{\qu\varphi}{\Pcal}{\Pcal'} extending $\varphi$:
   \begin{itemize}
   \vspace{-0.7em}
   \item[(a)] \PL\ admits a linear topology such that $\dom\varphi$
   is linearly open. Moreover, there exists a line $l\in\Lcal$ such that
   $(\dom\varphi \cap\l)^\varphi$ contains non--empty open set with
   respect to some linear topology of \PLL.
   \vspace{-0.7em}
   \item[(b)] \PL\ and \PLL\ are spaces of the same finite order $N$.
   Moreover, $x\cap\dom\varphi\neq\emptyset$ implies
   $\#(x\cap\dom\varphi) > \frac{2}{3}N \mbox{ for all }x\in\Lcal$.
   \end{itemize}
\end{theo}
Theorem \ref{frank} generalizes a result of {\sc Lenz \cite{lenz57b}} for
real projective spaces. There are also several theorems in the literature
giving local characterizations of embeddings of projective planes and spaces.
Refer to {\sc Bogn\'ar} and {\sc Kert\'esz} \cite{bogn-kert86}, {\sc
Brezuleanu} and {\sc R\u{a}dulescu} \cite{brez-radu85}, and the references
given in \cite{frank92}.

\section{Other Spaces}

\subsection{Partial Linear Spaces and Linear Spaces}

The geometric conditions (L1), (L2), (L1') etc. still do make sense in a
wider context.

Let \PL\ be a pair consisting of a {\em point set}\/ $\Pcal$ and a set
$\Lcal\subset 2^\Pcal$ whose elements are called {\em lines}. If any two
distinct points are on at most one common line, and if each line contains at
least two points, then \PL\ is called a {\em partial linear space}. A {\em
linear space}\/ is a partial linear space such that any two distinct points
are on a common line.

{\sc S\"orensen} \cite{soer85} has discussed mappings between linear spaces
\PL\ and \PLL\ satisfying (L1). Given $\Mcal\subset\Pcal$ the {\em trace
space}\/ $(\Mcal,\Lcal_\Mcal)$, which is defined by putting
   \begin{displaymath}
   \Lcal_\Mcal:=\{l\cap\Mcal\mid l\in\Lcal,\,\#(l\cap\Mcal)\geq 2\},
   \end{displaymath}
is also a linear space. Thus, if a linear mapping is not globally defined,
one can go over to the trace space on its domain. Hence it means no loss
of generality, that all mappings in \cite{soer85} are assumed to be globally
defined.

Examples of semicollineations of linear spaces and necessary conditions for a
semicollineation of linear spaces to be a collineation are due to {\sc
Kreuzer} \cite{kreu96}. There is a widespread literature on embeddings of
linear spaces in projective spaces \cite[Ch.\ 6]{buek95}.

\subsection{Affine Spaces}

The following seems to be part of the folklore: Suppose that a globally
defined mapping $\varphi$ of the affine space on $\Vbf$ (over $K$ with $\#
K>2$) in the affine space on $\Vbf'$ is satisfying (L1). If $\im\varphi$
contains a triangle, then $\varphi-0^\varphi$ is a semilinear mapping
$\Vbf\rightarrow\Vbf'$.

This result follows from Theorem \ref{frank} by going over to the projective
closures and, for infinite order, by applying their cofinite topology. See
{\sc S\"orensen} \cite{soer85} for a direct proof.

Globally defined mappings of affine spaces with property (WL1) have been
discussed by {\sc Zick} \cite{zick81a}, \cite{zick81b}. Since such a mapping
needs not preserve parallelity of lines, {\em fractional weak semilinear
mappings}\/ have been introduced to obtain an algebraic description.

Even an {\em embedding} $\varphi$ of the affine space on $\Vbf$ in the affine
space on $\Vbf'$ has in general not the property that $\varphi-0^\varphi$ is
a weak semilinear mapping. One class of examples is given by spaces over
\GF2. A more interesting example is the embedding of an affine plane over
$\mbox{GF}(3)$ in the complex affine plane which is based on the nine
inflections of an elliptic cubic curve, another example is due to {\sc Benz}
\cite[p.\ 113]{benz92}. See {\sc Bichara} and {\sc Korchm\'{a}ros}
\cite{bich-korc80}, {\sc Carter} and {\sc Vogt} \cite{cart-vogt80a},
\cite{cart-vogt80b}, {\sc Limbos} \cite{limb80a}, \cite{limb80b},
\cite{limb81}, \cite{limb82}, {\sc Ostrom} and {\sc Sherk}
\cite{ostr-sher64}, {\sc Rigby} \cite{rigb65}, {\sc Thas} \cite{thas70}. {\sc
Schaeffer} \cite{scha80} has given a local characterization of some
embeddings which can also be found in \cite[3.2]{benz92}.

\subsection{Grassmann Spaces}

Write $^d\Pcal$ for the set of all $d$--flats (i.e.\ $d$--dimensional
subspaces) of an $n$--dimensional projective space \PL\ $(1\leq d\leq n-2)$
and $^d\Lcal$ for the set of all pencils of $d$--flats. Then
$({}^d\Pcal,{}^d\Lcal)$ is a partial linear space, called {\em Grassmann
space}. Two distinct ``points'' $^dX,{}^dY\in {}^d\Pcal$ are ``collinear''
if, and only if, $^dX\cap{}^dY$ is a $(d-1)$--flat. The Grassmann space
$(^d\Pcal,{}^d\Lcal)$ is covered by one--dimensional projective spaces,
namely the pencils of $d$--flats.

A mapping $\pabb{\varphi}{{}^d\Pcal}{\Pcal'}$ in the point set of a
projective space \PLL\ is called {\em linear}, if the restriction of
$\varphi$ to each pencil is a linear mapping of projective spaces
\cite{havl-I81}.

If \PL\ is pappian, then each linear mapping
$\pabb{\varphi}{{}^d\Pcal}{\Pcal'}$ determines a mapping
$\pabb{\widehat\varphi}{G_{n,d}}{\Pcal'}$ of the associated Grassmann variety
$G_{n,d}$ (cf.\ {\sc Burau} \cite{bura61}) such that the restriction to each
line $l\subset G_{n,d}$ is linear. By {\sc Havlicek} \cite{havl-II81}, this
$\widehat\varphi$ extends to a unique linear mapping of the ambient space of
$G_{n,d}$. A weaker result is due to {\sc Wells} \cite{well83}.

If \PL\ is not pappian, then Grassmann varieties are not available due to a
non--existence theorem \cite{havl-II81}, and little seems to be known here.

A {\em geometric hyperplane}\/ of $({}^d\Pcal,{}^d\Lcal)$ is a subset of
${}^d\Pcal$ which intersects each pencil of $d$--flats in exactly one or in
all elements. The geometric hyperplanes are exactly the exceptional sets of
linear mappings $\pabb{\varphi}{{}^d\Pcal}{\Pcal'}$ with $\#\im\varphi=1$.
Hence in the pappian case all geometric hyperplanes arise as hyperplane
sections of the associated Grassmann variety \cite[4.5]{havl-II81}. All
geometric hyperplanes in the non--pappian case have been descibed by {\sc
Hall} and {\sc Shult} \cite{hall-shul93}.

It seems that {\sc Eckhart} \cite{eckh26} and {\sc Rehbock} \cite{rehb26a},
\cite{rehb26b} have been the first geometers to discuss linear mappings of
the Grassmann space formed by the lines of the real projective $3$--space.
See {\sc Brauner} \cite{brau73}, \cite{brau75}, \cite{brau83}, {\sc Havlicek}
\cite{havl83}, and {\sc L\"ubbert} \cite{lueb79} for further references
and {\em kinematic line mappings}. {\sc Zanella \cite{zane95}} has
investigated linear mappings of Grassmann spaces which are globally defined
and injective. He has given sufficient conditions for the image of such a
mapping to be projectively equivalent to the corresponding Grassmann variety
and rather sophisticated examples where this is not the case. Cf. also {\sc
Bichara} and {\sc Zanella} \cite{bich-zane97x}.

\subsection{Product Spaces}

Another partial linear space is the {\em product space}\/ of two projective
spaces. If both spaces have isomorphic commutative ground fields, then their
product space can be represented as a {\em Segre variety}. Linear mappings of
a product space in a projective space may be defined as for Grassmann spaces.
A thorough discussion of globally defined injective linear mappings can be
found in a paper by {\sc Zanella} \cite{zane96}. Unfortunately, we cannot
give the details here due to the lack of space. Let us just remark that the
situation is more complicated than for Grassmann spaces. Cf. also {\sc
Bichara} and {\sc Zanella} \cite{bich-zane97x}.

\subsection{Lattice Geometries}

A significant generalization of linear mappings to {\em projective}\/ and
{\em affine lattice geometries} has been given by {\sc Pfeiffer}
\cite{pfei93}, \cite{pfei97}, and {\sc Pfeiffer} and {\sc Schmidt}
\cite{pfei-schm95}. See furthermore {\sc Faigle} \cite{faig80}.

\abstand
\noindent
Abteilung f\"ur Lineare Algebra und Geometrie, Technische Universit\"at,
Wiedner Hauptstra{\ss}e 8--10, A--1040 Wien, Austria.\\
EMAIL: {\tt havlicek@geometrie.tuwien.ac.at}
\end{document}